\documentclass[letterpaper,11pt,oneside]{amsart}
\usepackage{amsmath}
\usepackage{amssymb}
\usepackage{latexsym}
\usepackage{amsfonts}
\usepackage[english,activeacute]{babel}
\usepackage{multicol}
\usepackage{oldgerm}
\usepackage{enumerate}
\usepackage{hyperref}
\usepackage{mathtools}
\usepackage{tikz}
\usepackage{verbatim}
\usepackage{xcolor} 
\usepackage[mathscr]{euscript}

\usepackage{comment}

\usepackage{float}

\newtheorem{teo}{Theorem}
\newtheorem{lema}{Lemma}

\newcommand{\oeis}[1]{\href{http://oeis.org/#1}{\underline{#1}}}

\def\BigConstant#1#2{#1 \times 10^{#2}}    
\def\F#1#2{F^{(#1)}_{#2}}    
\def\A{\alpha}               
\def\g{\gamma}               

\def\K{\mathbb{K}}

\def\Q{\mathbb{Q}}

\begin{document}

\title[$\F{k}{n}\F{k}{m} = P_{\ell}$]{On Pell numbers representable as product of two generalized Fibonacci numbers}

\author{Jhon J. Bravo}
\address{Departamento de Matem\'aticas\\ Universidad del Cauca\\ Calle 5 No 4--70\\Popay\'an, Colombia.}
\email{jbravo@unicauca.edu.co}

\author{Pranabesh Das}
\address{Department of Mathematics\\ Xavier University of Louisiana\\ 1 Drexel Drive, New Orleans, Louisiana, 70125}
\email{pranabesh.math@gmail.com, pdas@xula.edu}

\author{Jose L. Herrera}
\address{Departamento de Matem\'aticas\\ Universidad del Cauca\\  Calle 5 No 4--70\\Popay\'an, Colombia.}
\email{joseherrera@unicauca.edu.co}

\author{John C. Saunders}
\address{Department of Mathematical Sciences\\ Middle Tennessee State University\\ 1301 E. Main Street, Murfreesboro, Tennessee, 37130}
\email{John.Saunders@mtsu.edu}
	
\date{\today}
	
\begin{abstract}
A generalization of the well-known Fibonacci sequence is the $k$-Fibonacci sequence with some fixed integer $k\ge 2$. The first $k$ terms of this sequence are $0,0, \ldots, 1$, and each term afterwards is the sum of the preceding $k$ terms. In this paper, we find all Pell numbers that can be written as product of two $k$-Fibonacci numbers. The proof of our main theorem uses lower bounds for linear forms in logarithms, properties of continued fractions, and a variation of a result of Dujella and Peth\H{o} in Diophantine approximation. This work generalizes a prior result of Alekseyev which dealt with determining the intersection of the Fibonacci and Pell sequences, a work of Ddamulira, Luca and Rakotomalala who searched for Pell numbers which are products of two Fibonacci numbers, and a result of Bravo, G\'omez and Herrera who found all Pell numbers appearing in the $k$-Fibonacci sequence.    


\medskip


		
\end{abstract}
	
\maketitle

\section{Introduction}\label{sec:Introduction}

\noindent For $k \geq 2$, we focus on the $k$-generalized Fibonacci sequence or, for simplicity, the $k$-Fibonacci sequence $\F{k}{} \coloneqq (\F{k}{n})_{n\ge 2-k}$ which is a linear recurrence sequence of order $k$ defined by the relation  
\[
\F{k}{n}=\F{k}{n-1}+\F{k}{n-2}+\cdots+\F{k}{n-k}\quad \text{for all}\quad n\geq 2,
\]
with $\F{k}{-(k-2)}=\F{k}{-(k-3)}=\cdots=\F{k}{0}=0$ and $\F{k}{1}=1$ as initial conditions. We usually call $\F{k}{n}$ the $n$th $k$-Fibonacci number.  The name of this family of sequences comes from the fact that $\F{k}{}$ coincides with the famous Fibonacci sequence $(F_n)_{n\geq 0}$ when $k=2$. A compendium of the most important and interesting properties of the Fibonacci numbers can be found, for example, in \cite[\oeis{A000045}]{Sloane} and Koshy's book \cite{Kos01}. The $k$-Fibonacci sequence for the particular cases $k=3,4$ can also be found in \cite{Sloane} as the sequences \oeis{A000073} and \oeis{A000078}, respectively.

Regarding the sequence $\F{k}{}$, it is well-known that its first $k+1$ non-zero terms are powers of two; namely
\begin{equation}\label{primeros-k}
\F{k}{n} = 2^{\max\{0,n-2\}} \quad\text{for all}\quad  1 \le n \le k+1.
\end{equation}
A wide range of recent studies have focused on the problem of determining all $k$-Fibonacci numbers that admit a special representation. For instance, Bravo and G\'omez \cite{BG2016} determined all terms in $F^{(k)}$ that are one less than a power of $2$. In another work, Bravo, G\'omez and Luca \cite{BGL1} identified all terms in $F^{(k)}$  that can be written as the sum of two repdigits, while Bravo and Luca \cite{BravoLuca2015} characterized all repdigits that can be expressed as the sum of of two $k$-Fibonacci numbers. In addition, Bravo, G\'omez and Luca \cite{BGL2016} found all powers of 2 that are the sum of two $k$-Fibonacci numbers. 

Another important sequence, closely related to the Fibonacci sequence, is the Pell sequence. The Pell sequence $(P_{\ell})_{\ell\geq 0}$ is defined by the recurrence 
\[
P_{\ell}=2P_{\ell-1}+P_{\ell-2}  \quad \text{for all} \quad \ell\geq 2,
\]
with $P_0=0$ and $P_1=1$ as initial conditions. This sequence corresponds to entry \oeis{A000129} in \cite{Sloane}. For the beauty and rich applications
of these numbers and their relatives one can see Koshy's book \cite{Kos14}.

Much research has been devoted to Diophantine equations involving Pell and $k$-Fibonacci numbers. For example, Alekseyev \cite{A1} proved that the only integers belonging to both the Fibonacci and Pell sequences are $0,1,2,5$, and also established analogous results for other Lucas sequences. Bravo, G\'omez and Herrera \cite{BGH-19} extended these findings to the Pell and $k$-Fibonacci sequences, showing that for every $k\geq 2$, the only nontrivial common terms are $5$ and $29$. Furthermore, Bravo, Herrera and Luca \cite{BHL_2021_JNT} provided a broader generalization by determining all integers that appear in both generalized Pell and generalized Fibonacci sequences.

This line of research, originally focused on finding intersections of two linear recurrence sequences, has been extended to the broader problem of determining all terms in a given sequence that can be expressed as products of terms from another sequence. In particular, such results have been given on the $k$-Fibonacci, Pell, and other linear recurrence sequences. For example, in 2016, Ddamulira, Luca and Rakotomalala \cite{DdamuLucaRako2016} proved that the only Fibonacci numbers that are the product of two Pell numbers are $1,2,5,144$, and the only Pell numbers that are the product of two Fibonacci numbers are $1,2,5,169$. Likewise, Da\c{s}demir and Varol \cite{Dasdemir} found all Pell numbers that are the product of two Lucas numbers, and all Lucas numbers that are the product of two Pell numbers. Also, Rihane, Akrour and Habibi \cite{Rihane1} found all Fibonacci numbers that are the product of three Pell numbers and all Pell numbers that are the product of three Fibonacci numbers. Moreover, Rihane \cite{Rihane2} proved the following theorem, giving all the terms in the $k$-Fibonacci sequences that are products of two terms in the Pell sequence.

\begin{teo}[Rihane \cite{Rihane2}]
All positive integer solutions $(\ell,m,k,n)$ of the Diophantine equation
\begin{equation}
P_{\ell}P_{m}=\F{k}{n},
\end{equation}
with $k\geq 3$ and $2\leq \ell \leq m$, are
\[
(\ell,m,k,n)=(2,2,k,4),(2,4,3,7),(4,7,8,13),(3,8,9,13).
\]
\end{teo}

In this paper, we pursue the opposite problem and search for Pell numbers that can be expressed as the product of two $k$-Fibonacci numbers. In this sense, our work may be viewed as an extension of the results established in \cite{A1,BGH-19,DdamuLucaRako2016}. 

Our main result is as follows.

\begin{teo}\label{Main_Theorem}
All positive integer solutions $(n,m,k,\ell)$ of the Diophantine equation
\begin{equation}\label{eq:Main_Equation}
\F{k}{n}\F{k}{m} = P_{\ell},
\end{equation}
with $k\geq 2$ and $n\geq m\geq 3$, are
\[
(n,m,k,\ell)=(7,7,2,7),(6,6,3,7),(15,3,5,12).
\]
\end{teo}
In the previous theorem, we can establish the condition $m\geq 3$ because if $m=1$ or $m=2$, the original  equation reduces to $\F{k}{n}= P_{\ell}$, and this case was already solved by Alekseyev \cite{A1} and by Bravo, G\'omez and Herrera \cite{BGH-19}, as previously mentioned.


\section{Preliminary Results}\label{sec:Preliminary_Results}

\noindent In this section, we present the notation and background required for the proof of Theorem \ref{Main_Theorem}. In particular, we review some results concerning Pell and generalized Fibonacci numbers, along with the mathematical tools that will be employed throughout the paper. 

\subsection{The Pell sequence}
An explicit Binet formula for the Pell sequence is well known. Namely, for all $\ell\geq 0$, we have that
\begin{equation}\label{Binet.Pell}
P_{\ell}=\frac{\gamma^{\ell}-\overline{\gamma}^{\ell}}{\gamma-\overline{\gamma}}=\frac{\gamma^{\ell}-\overline{\gamma}^{\ell}}{2\sqrt{2}},
\end{equation}
where $(\gamma,\overline{\gamma}):=(1+\sqrt{2}, 1-\sqrt{2})$ are the roots of the characteristic equation $x^{2}-2x-1=0$. In particular, \eqref{Binet.Pell} easily implies that the inequality
\begin{equation}\label{desPell}
  \gamma^{\ell-2} \leq P_{\ell} \leq \gamma^{\ell-1} \quad \text{holds for all} \quad \ell\geq 1. 
\end{equation}
Additionally,  in view of \eqref{Binet.Pell}, we can write
\begin{equation}\label{AproximacionPell}
P_{\ell} = \frac{\gamma^{\ell}}{2\sqrt{2}}+\xi(\ell) \quad \text{where} \quad |\xi(\ell)|<1/5 \quad \text{for all} \quad \ell\geq 1. 
\end{equation}
There are several well-known results concerning the Pell sequence. For example, a result of Ljunggren \cite{Ljunggren42} shows that $P_{\ell}$ is a perfect square only when $\ell=0,1$ or $7$. Later, Peth\H{o} \cite{Petho91} proved that these are, in fact, the only perfect powers in the entire Pell sequence (see also \cite{Cohn96}). We state their result which will be used in the sequel.

\begin{lema}\label{pethocohn}
The Diophantine equation $P_{\ell}=x^m$ in positive integers $\ell$, $x$ and $m$ with $m\geq 2$ has only the solutions $(\ell, x, m)=(1, 1, m)$ and $(7,13,2)$.
\end{lema}


\subsection{On $k$-Fibonacci numbers}
We begin by mentioning that the $k$-Fibonacci sequence is a linear recurrence sequence with characteristic polynomial
\[
\Psi_k(z) \coloneqq z^k-z^{k-1}-\cdots-z-1,
\]
which is irreducible over $\Q$. This polynomial has exactly one positive real root $\A:=\A(k)$ outside the unit complex disk (see \cite{Mi,Mil}), while the other roots are strictly inside the unit circle. More precisely, $2(1-2^{-k})<\A(k)<2$ (see \cite[Lemma 2.3]{HY} or \cite[Lemma 3.6]{DAW}). The root $\A(k)$ is known as the \emph{dominant root} of $\F{k}{}$ and in what follows we shall omit the dependence on $k$ of $\A$. 

We now present a lemma that summarizes the main properties of the $k$-Fibonacci sequence (see \cite{BGL2016, BL1, DD}). This result will be used repeatedly in the proof of Theorem~\ref{Main_Theorem}.

\begin{lema}[Properties of $\F{k}{}$]\label{properties-F}
    Let $k\geq 2$ be an integer. Then 
    \begin{itemize} 
        \item[$(a)$]  $\alpha^{r-2} \leq \F{k}{r}\leq \alpha^{r-1}$ for all $r\geq 1$.
        
        \item[$(b)$] $\F{k}{}$ satisfies the following ``Binet--like" formula
        \[
        \F{k}{r}=\sum_{i=1}^{k}f_{k}(\alpha_i)\alpha_i^{r-1},
        \]     
        where $\alpha=\alpha_1,\ldots,\alpha_{k}$ are the roots of $\Psi_{k}(z)$ and 
        \[
        f_{k}(z) \coloneqq \frac{z-1}{2+(k+1)(z-2)}.
        \]
        
        \item[$(c)$] $1/2< f_{k}(\alpha) <3/4$ and $| f_{k}(\alpha_i)|< 1$ for $2\leq  i \leq k$.
        
        \item[$(d)$] $|\F{k}{r}-f_{k}(\alpha)\alpha^{r-1}|<1/2$  holds for all $r\geq 2-k$.
        
    \end{itemize}
\end{lema}

The following estimate due to Bravo, G\'omez and Luca \cite{BGL1} is currently one of the key points in addressing the large values of $k$ when solving Diophantine equations involving terms of generalized Fibonacci numbers (see also \cite{BGH-19}).
	
\begin{lema}\label{E-large-k}
Let $k\geq 2$ and suppose that $r<2^{k/2}$. Then
\[
F_{r}^{(k)}=2^{r-2}(1+\zeta_{r})\quad \text{where}\quad |\zeta_{r}|<\frac{2}{2^{k/2}}.
\]
\end{lema}

\subsection{Linear forms in logarithms}

To prove Theorem \ref{Main_Theorem}, we use lower bounds for linear forms in logarithms of algebraic numbers. In this subsection, we present the main tool for this purpose. Before that, we recall some important concepts from algebraic number theory.

Let $\eta$ be an algebraic number of degree $d$, and let $\sum_{0\le j\le d}a_{j}x^{j}$ be its minimal primitive polynomial in $\mathbb{Z}[x]$. The logarithmic height $h(\eta)$ of $\eta$ is given by
$$
h(\eta)=\dfrac{1}{d}\left(\log(|a_{d}|)+\sum_{i=1}^{d}\log\max\bigl\{|\eta_{i}|,1\bigr\}\right),
$$
where $\eta_1, \eta_2, \ldots, \eta_d$ are the conjugates of $\eta$. In particular, if $\eta=p/q \in \mathbb{Q}$ is in the lowest terms with $q\ge 1$, then $h(\eta)=\log \max \bigl\{|p|,q\bigr\}$. The following properties \cite[Property 3.3]{Waldshmidt} of the logarithmic height will be used with or without further reference as and when needed.

\begin{itemize}
		\item $h(\eta_{1}\pm \eta_{2})\leq h(\eta_{1})+h(\eta_{2})+\log 2$;
		\item $h(\eta_{1}\eta_{2}^{\pm 1})\leq h(\eta_{1})+h(\eta_{2})$;
		\item $h(\eta^{r})=|r|h(\eta)$, $r\in \mathbb{Z}$;
\end{itemize}

For later estimates, it will be important to take into account that
\begin{equation}\label{Estimation_h(Alpha)_h(f(Alpha))}
    h(\A) = (\log \A) / k < (\log 2)/k 
    \qquad \text{and} \qquad
    h(f_{k}(\A))< 2\log k \quad\text{for all}\quad k \ge 2. 
\end{equation}
See \cite{BG2016} for further details of the proof of \eqref{Estimation_h(Alpha)_h(f(Alpha))}.

Let $\eta_{1},\ldots,\eta_{t}$ be nonzero elements of a number field $\mathbb{K}$ of degree $D$, and let $b_{1},\ldots,b_{t}$ be integers. Set
	$$
	B\geq \max \bigl\{|b_1|,\ldots,|b_{t}|\bigr\} \quad {\rm and} \quad \Lambda =\eta_{1}^{b_1}\cdots \eta_{t}^{b_t}-1.
	$$
Let $A_{1},\ldots,A_{t}$ be real numbers such that
	$$
	A_{i}\geq \max \bigl\{Dh(\eta_{i}),\left|\log \eta_{i}\right|,0.16\bigr\}\quad\text{for all}\quad i=1,\ldots,t.
	$$
With the previous notation, the main result of Matveev~\cite{Matveev} implies the following estimate.
\begin{teo}\label{Matveev}
Assume that $\Lambda$ is nonzero. If $\mathbb{K}$ is real, then
$$
|\Lambda| > \exp \left(  -1.4\cdot 30^{t+3}\cdot t^{4.5}\cdot D^{2}\left(1+\log D\right)\left(1+\log B\right)A_{1}\cdots A_{t} \right).
$$
\end{teo}

\subsection{Reduction tools}

Theorem \ref{Matveev} is a fundamental tool for solving certain exponential Diophantine equations. However, the bounds obtained from Matveev’s result are typically quite large, making them unsuitable for subsequent computational searches. To reduce these bounds, we employ a result from the theory of continued fractions. The following lemma is a slight variation of a result originally due to Dujella and Peth\H{o}~\cite{DP}. Specifically, we use the version presented by Bravo, G\'omez and Luca (see \cite[Lemma 1]{BGL2016}).

\begin{lema}\label{reduce}
Let $M$ be a positive integer and let $A, B, \mu$ and $\tau$ be given real numbers with $A>0$ and $B>1$. Assume that $p/q$ is a convergent of the continued fraction of $\tau$ such that $q>6M$. Let $\epsilon \coloneqq \left\|\mu q\right\|-M\left\|\tau q\right\|$, where $\left\|\cdot\right\|$ denotes the distance from the nearest integer. If $\epsilon>0$, then there is no solution to the inequality
$$
0<|u\tau-v+\mu|<AB^{-w}
$$
in positive integers $u$, $v$ and $w$ with $u\leq M$ and $w\geq \log (Aq/\epsilon)/\log B$.
\end{lema}

We finish this subsection with an analytical tool and a simple fact concerning the exponential function (for a proof, see \cite[Lemma 7]{GL} and \cite[Lemma 8]{BHL_2021_JNT}, respectively). We list them as lemmas for further reference.  

\begin{lema}\label{GuzmanLemma}
If $s\geq 1$ is an integer and $x$ and $T$ are real numbers such that
\[
T > (4s^2)^s \quad\text{and}\quad \frac{x}{\log^{s}x} < T, \qquad \text{then} \qquad x < 2^s T \log^{s} T.
\]
\end{lema}

\begin{lema}\label{exponential}
For any nonzero real number $z$, we have
\begin{enumerate}
\item[$(a)$] $0< z < |e^z - 1|$.
\item[$(b)$] If $z<0$ and $|e^z - 1| <1/2$, then $|z| < 2\, |e^z - 1|$.
\end{enumerate}
\end{lema}

\section{Proof of Theorem \ref{Main_Theorem}}

\noindent Assume throughout that $(n,m,k,\ell)$ is a solution of the equation \eqref{eq:Main_Equation}. First, note that if $k=2$, the resulting equation was completely solved by Ddamulira, Luca and Rakotomalala in \cite{DdamuLucaRako2016}, where it was shown that $(n,m,k,\ell)=(7,7,2,7)$ is the only solution with $m\geq 3$. Therefore, from now on, we may assume that $k\geq 3$. On the other hand, if $n=m$, then equation \eqref{eq:Main_Equation} becomes $(F_{n}^{(k)})^2=P_{\ell}$. Since $n\geq 3$, Lemma \ref{pethocohn} implies that $\ell =7$, and consequently $F_{n}^{(k)}=13$. This yields $k=3$ and $n=6$. Hence, the only solution to equation \eqref{eq:Main_Equation} with $k\geq 3$ and $n=m$ is $(n,m,k,\ell)=(6,6,3,7)$. Thus, in the remainder of the proof, we may assume that $n>m\geq 3$.

Now, if $n\leq k+1$, then $F_{n}^{(k)}=2^{n-2}$ and $F_{m}^{(k)}=2^{m-2}$ by \eqref{primeros-k} and so equation \eqref{eq:Main_Equation} becomes $2^{n+m-4}=P_{\ell}$, which can be seen to have no solutions in view of Lemma \ref{pethocohn}. Thus, for the rest of the paper, we assume that $n\geq k+2$ which gives $n\geq 5$ and $n+m\geq 8$. It then follows from the above and \eqref{eq:Main_Equation} that  $P_{\ell} \geq F_{5}^{(3)}F_{3}^{(3)}=7\cdot 2=14$ implying $\ell \geq 5$. 

We now want to establish a relationship between the variables of \eqref{eq:Main_Equation}. For this purpose, we use \eqref{eq:Main_Equation}, \eqref{desPell} and Lemma \ref{properties-F}$(a)$ to get 
\[
\alpha^{n+m-4}\leq F_{n}^{(k)}F_{m}^{(k)}=P_{\ell}< \gamma^{\ell-1}  
\]
and
\[
\gamma^{\ell-2}<P_{\ell}=F_{n}^{(k)}F_{m}^{(k)}\leq \alpha^{n+m-2}.  
\]
Thus
\[
(n+m-4)\left(\frac{\log\alpha}{\log\gamma}\right)+1<\ell<(n+m-2)\left(\frac{\log\alpha}{\log\gamma}\right)+2.
\]
Taking into account that $n+m\geq 8$ and using the fact that $1.89<\alpha(3)\leq \alpha<2$, we conclude that
\begin{equation}\label{des-var}
\frac{2(n+m)}{5}<\ell<\frac{9(n+m)}{10},  
\end{equation}
which is an estimate on $\ell$ in terms of $n+m$.

On the other hand, by Lemma \ref{properties-F}$(d)$, we can write 
\begin{equation}\label{error_Fr}
\F{k}{r}=f_{k}(\alpha)\alpha^{r-1}+e_{k}(r)
\qquad \text{where} \qquad 
|e_{k}(r)|<1/2, 
\end{equation}
and a simple computation gives
\begin{equation}\label{aux_error}
\F{k}{r}=f_{k}(\alpha)\alpha^{r-1} ( 1 +  X_{r}) 
\qquad \text{where} \qquad 
|X_{r}| < \frac{2}{\alpha^{r}}.
\end{equation}
At this point, we use \eqref{AproximacionPell} and \eqref{aux_error} to rewrite \eqref{eq:Main_Equation} as
\[
f_{k}^{2}(\alpha)\alpha^{n+m-2}-\frac{\g^{\ell}}{2\sqrt{2}}=\xi(\ell)-f_{k}^{2}(\alpha)\alpha^{n+m-2}X,
\]
where $X=X_{n}+X_{m}+X_{n}X_{m}$ for which $|X|<8/\alpha^{m}$. Taking absolute values and dividing both sides of the above relation by $f_{k}^{2}(\alpha)\alpha^{n+m-2}$, we obtain the inequality
\begin{equation}\label{forma1}
|\Lambda_{1}|:=\left|\g^{\ell} \cdot \alpha^{-(n+m-2)} \cdot (2\sqrt{2}f_{k}^{2}(\alpha))^{-1}-1\right|<\frac{10}{\alpha^{m}}.    
\end{equation}
In the first application of Theorem \ref{Matveev}, we take $t:=3$ and the parameters
\[
(\eta_{1},b_{1}):=(\g,\ell), \quad (\eta_{2},b_{2}): =(\alpha,-(n+m-2)) \quad \text{and} \quad (\eta_{3},b_{3}):=(2\sqrt{2}f_{k}^{2}(\alpha),-1).
\]
Let us start by noticing that $\eta_{1},\eta_{2}$ and $\eta_{3}$ are real numbers belonging to $\mathbb{K}:=\mathbb{Q}(\alpha,\sqrt{2})$ for which $D:=[\mathbb{K}:\mathbb{Q}]\leq 2k$. The fact that $\Lambda_{1}\neq 0$ is justified as follows. Suppose, for the sake of contradiction, that $\Lambda_{1}=0$. Then
\begin{equation}\label{absurdo1}
f_{k}^{2}(\alpha)\alpha^{n+m-2}=\frac{\g^{\ell}}{2\sqrt{2}}.
\end{equation}
To see that \eqref{absurdo1} is not possible, we argue as follows. Let $\mathbb{L}=\Q(\A_{1},\ldots,\A_{k},\g)$ be the normal closure of $\K$ and let $\sigma_{1},\ldots,\sigma_{k}$ be elements of $\text{Gal}(\mathbb{L}/\Q)$ such that $\sigma_{i}(\A)=\A_{i}$. Since $k\geq 3$ and $\g$ has two conjugates, the elements $\sigma_{1}(\g),\ldots,\sigma_{k}(\g)$ cannot all be distinct. Therefore, there exist indices $i\neq j$ in $\left\{1,2,\ldots,k\right\}$ such that $\sigma_{i}(\g)=\sigma_{j}(\g)$. Note that, if we put $\sigma_{j}^{-1}(\A_{i})=\A_{s}$, then $s\neq 1$ because $\sigma_{j}(\A_{1})=\A_{j}\neq\A_{i}$. Applying $\sigma_{j}^{-1}\sigma_{i}$ to the relation \eqref{absurdo1}, we get 
\[
f_{k}^{2}(\A_{s})\A_{s}^{n+m-2}=\frac{\g^{\ell}}{2\sqrt{2}}.
\]
Taking absolute value in the above relation we see that this is not possible since its right-hand side exceeds $1$ for all $\ell\geq 5$, while its left-hand side is smaller than $1$. Thus, $\Lambda_{1}\neq 0$.

On the other hand, $h(\eta_{1})=(\log\g)/2$ and $h(\eta_{2})=(\log\alpha)/k<(\log 2)/k$, so one may take $A_{1}:=k\log\g$ and $A_{2}:=2\log 2$.  To estimate $h(\eta_{3})$, we use the properties of the logarithmic height to get
\[
h(\eta_{3})\leq h(2\sqrt{2})+2h(f_{k}(\alpha))\leq \frac{3\log 2}{2}+4\log k\leq 5\log k
\]
for all $k\geq 3$, where we used the fact that $h(f_{k}(\alpha))<2\log k$ as mentioned before. Hence, we can choose $A_{3}:=10k\log k$.  Finally, since $\ell<2n$ by \eqref{des-var}, it is appropriate to set $B:=2n$. Then, Theorem \ref{Matveev} tells us that
\begin{equation}\label{E1}
\log |\Lambda_{1}|>-3.83\times 10^{13}k^{4}\log^{2}k\log n,
\end{equation}
where the inequalities $1+\log(2k)\leq 2.6 \log k$ for all $k\geq 3$ and $1+\log(2n)\leq 2.1 \log n$ for all $n\geq 5$ have been applied. Taking logarithms in \eqref{forma1} and comparing the resulting inequality with \eqref{E1}, it follows that
\begin{equation}\label{E2}
m\log\alpha<3.9\times 10^{13}k^{4}\log^{2}k\log n.
\end{equation}
Let us now get a second linear form in real logarithms. To this end, we use \eqref{Binet.Pell} and \eqref{aux_error} again to rewrite \eqref{eq:Main_Equation} this time in the form
\[
f_{k}(\alpha)\alpha^{n-1}F_{m}^{(k)}-\frac{\g^{\ell}}{2\sqrt{2}}=\xi(\ell)-f_{k}(\alpha)\alpha^{n-1}F_{m}^{(k)}X_{n},
\]
where $|X_{n}|<2/\alpha^{n}$. Taking absolute values on both sides of the above relation and dividing it across by $f_{k}(\alpha)\alpha^{n-1}F_{m}^{(k)}$, we obtain
\begin{equation}\label{forma2}
|\Lambda_{2}|:=\left|\g^{\ell} \cdot \alpha^{-(n-1)} \cdot (2\sqrt{2}f_{k}(\alpha)F_{m}^{(k)})^{-1}-1\right|<\frac{3}{\alpha^{n}},
\end{equation}
where we used that $F_{m}^{(k)}\geq F_{3}^{(3)}=2$. In the second application of Theorem \ref{Matveev}, we take $t:=3$ and 
 \[
(\eta_{1},b_{1}):=(\g,\ell), \quad (\eta_{2},b_{2}): =(\alpha,-(n-1)) \quad \text{and} \quad (\eta_{3},b_{3}):=(2\sqrt{2}f_{k}(\alpha) F_{m}^{(k)},-1).
\]
As before $\mathbb{K}:=\mathbb{Q}(\alpha,\sqrt{2})$ contains $\eta_{1},\eta_{2}$ and $\eta_{3}$ and has $D:=[\mathbb{K}:\mathbb{Q}]\leq 2k$. The choices of $A_1$, $A_2$ and $B$ are the same as in the first application. To see why $\Lambda_{2}\neq 0$, note that
otherwise, we would get the relation 
\[
\g^{\ell}=2\sqrt{2}f_{k}(\alpha)\alpha^{n-1}F_{m}^{(k)}.
\]
Applying now the automorphism used in the first application of Matveev’s theorem, which fixes $\g$ and maps $\alpha$ to another conjugate, we obtain that for some $s\geq 2$,
\[
\g^{\ell}=2\sqrt{2}f_{k}(\alpha_{s})\alpha_{s}^{n-1}F_{m}^{(k)}.
\]
Taking absolute values on both sides yields
\[
\g^{\ell}<2\sqrt{2}F_{m}^{(k)}=2\sqrt{2}\frac{P_{\ell}}{F_{n}^{(k)}}<\frac{2\sqrt{2}}{7}\g^{\ell-1},
\]
where we used the fact that $F_{n}^{(k)}\geq F_{5}^{(3)}=7$. It then follows from the above that $\g<2\sqrt{2}/7$ which is not possible. Thus, $\Lambda_{2}\neq 0$.

Let us now estimate $h(\eta_{3})$. Applying the properties of $h(\cdot)$ and taking into account the inequality \eqref{E2}, we get
\begin{align*}
h(\eta_{3}) &\leq h(2\sqrt{2})+h(f_{k}(\alpha))+h(F_{m}^{(k)}) \\
& < \frac{3\log 2}{2}+2\log k+(m-1)\log \alpha\\
& < 4\times 10^{13}k^{4} \log^{2}k \log n. 
\end{align*}
Therefore, we can take $A_{3}:=8\times 10^{13}k^{5} \log^{2}k \log n$. Theorem \ref{Matveev} now tells us that
\begin{equation}\label{E3}
\log |\Lambda_{2}|>-3.1\times 10^{26}k^{8}\log^{3}k\log^{2}n,
\end{equation}
where the inequalities $1+\log(2k)\leq 2.6 \log k$ for all $k\geq 3$ and $1+\log(2n)\leq 2.1 \log n$ for all $n\geq 5$ have also been used. Comparing the above inequality \eqref{E3} with \eqref{forma2} implies that
\[
n<4.9\times 10^{26}k^{8}\log^{3}k\log^{2}n,
\]
and so
\begin{equation}\label{aplicGuzman}
\frac{n}{\log^{2}n.}<4.9\times 10^{26}k^{8}\log^{3}k.
\end{equation}
In order to get an upper bound for $n$ depending on $k$ we next use Lemma \ref{GuzmanLemma}. Indeed, taking $s:=2$, $x:=n$, $T=4.9\times 10^{26}k^{8}\log^{3}k$ and using the fact that 
\[
\log(4.9\times 10^{26}k^{8}\log^{3}k)<61.46+8\log k+3\log \log k\leq 65\log k
\]
for all $k\geq 3$, inequality \eqref{aplicGuzman} yields $n<8.3\times 10^{30}k^{8}\log^{5}k$. We record what we have proved so far as a lemma.
\begin{lema}\label{n-in-terms-ofk}
If $(n,m,k,\ell)$ is a solution of equation \eqref{eq:Main_Equation} with $k\geq 3$ and $n\geq k+2$, then the inequalities
\[
n<8.3\times 10^{30}k^{8}\log^{5}k \quad \text{and} \quad \ell<1.7\times 10^{31}k^{8}\log^{5}k 
\]
hold.
\end{lema}

\subsection{The case of small $k$}

Suppose that $k\in [3, 420]$. Here, we need to find better bounds for $n$, $m$ and $\ell$ than those implied by Lemma \ref{n-in-terms-ofk} for these values of $k$. To do this, we first let
\[
z_{1}:=\ell\log\g-(n+m-2)\log\alpha-\log(2\sqrt{2}f_{k}^{2}(\alpha))
\]
and observe that \eqref{forma1} can be rewritten as
\[
|e^{z_{1}}-1|<\frac{10}{\alpha^{m}}.
\]
Suppose that $m\geq 5$ and note that $z_{1} \neq 0$. If $z_{1} > 0$, we apply Lemma \ref{exponential}$(a)$ to deduce that $|z_1| < | e^{z_1} -1 | < 10/\alpha^{m}$.  If, on  the contrary, $z_1 < 0$, then $|e^{z_1} -1| <10/\alpha^{m} < 1/2$ since $m\ge 5$. Hence, from Lemma \ref{exponential}$(b)$ it follows that $|z_1| < 2| e^{z_1} -1 | <20/\alpha^{m}$. In either case, we obtain the bound $|z_1| < 20/\alpha^{m}$.  Replacing $z_1$ in the aforementioned inequality by its formula and dividing it across by $\log\alpha$, one can see that
\begin{equation}\label{DP:1}
0<\left|\ell\tau_k -(n+m-2)+ \mu_k \right| < A \cdot B^{-m},
\end{equation}
where we have put
\[
\tau_k:= \frac{\log\g}{\log\A},\quad \mu_k:=- \frac{\log(2\sqrt{2}f_k^{2}(\A))}{\log\A}, \quad A:= 32 \quad \text{and} \quad B:=\A.
\]
We also set $M_k:=\lfloor 1.7\times 10^{31} k^8 \log^5 k\rfloor$ which is an upper bound on $\ell$ by Lemma \ref{n-in-terms-ofk}. Observe that $\tau_k$ is irrational since $\g$ and $\A$ are multiplicatively independent. Applying now Lemma \ref{reduce} to inequality \eqref{DP:1} for each $k\in[3,420]$ we find, with help of \emph{Mathematica}, that $m \le 210$.

We now put
\[
z_2:= \ell\log\g - (n-1)\log\A -\log (2\sqrt{2}f_{k}(\A)F_{m}^{(k)}).
\] 
This allows us to rewrite inequality \eqref{forma2} as
\[
|e^{z_2} -1| < \frac{3}{\A^{n}}.
\]
Note that $z_2 \neq 0$. By analyzing the cases where $z_{2}$ is positive or negative, and taking into account that $|e^{z_2} -1| <3/\A^{n}< 1/2$ since $n\ge 5$, Lemma \ref{exponential}, along with the earlier argument, allows us to deduce that $|z_2| < 6/\A^{n}$. Replacing $z_2$ in the aforementioned inequality by its
formula and dividing it across by $\log\A$, we conclude that
\begin{equation}\label{DP:2}
0<\left|\ell \tau_k - (n-1) + \mu_{k,m} \right| <  A \cdot B^{-n},
\end{equation}
where we now put 
\[
\tau_k:= \frac{\log\g}{\log\A},\quad \mu_{k,m}:=- \frac{\log(2\sqrt{2}f_k(\A)F_{m}^{(k)})}{\log\A}, \quad A:= 10 \quad \text{and} \quad B:=\A.
\]
With this new choice and taking again $M_k:=\lfloor\BigConstant{1.7}{31} k^8 \log^5 k\rfloor$ as upper bound on $\ell$, we apply Lemma \ref{reduce} to inequality \eqref{DP:2} for each $k\in[3,420]$ and $m \in [3,210]$ obtaining that the possible solutions all have $n\leq 205$. This bound on $n$ allows us to conclude that $\ell <370$ by \eqref{des-var} and $k\leq 203$ since $k+2\leq n$. Finally, a brute force
search done in the range      
\[
3\le k \le 203,\quad
k + 2 \le n \le 205,
\quad 3\le m< n
\quad\text{and}\quad 
\quad 5\leq \ell< 370
\]
gives the sporadic solution $(n,m,k,\ell)=(15,3,5,12)$. 


\subsection{The case of large $k$}

Suppose now that $k>420$. In this case the following inequalities hold:
\[
m<n<8.3\times 10^{30}k^{8}\log^{5}k<2^{k/2}.
\]
At this point, we require the estimation from Lemma \ref{E-large-k} to prove that equation \eqref{eq:Main_Equation} has no
solutions for large values of $k$. Indeed, since $m<n<2^{k/2}$, Lemma \ref{E-large-k}, together with \eqref{eq:Main_Equation} and \eqref{AproximacionPell}, implies that
\[
2^{n+m-4}-\frac{\g^{\ell}}{2\sqrt{2}}=\xi(\ell)-2^{n+m-4}\zeta,
\]
where $\zeta=\zeta_{n}+\zeta_{m}+\zeta_{n}\zeta_{m}$ for which $|\zeta|<8/2^{k/2}$. Dividing both sides of the above equation by $2^{n+m-4}$ and taking absolute values we get 
\begin{equation}\label{forma3}
|\Lambda_{3}|:=\left|\g^{\ell}\cdot 2^{-(n+m-3)}\cdot (\sqrt{2})^{-1}-1\right |<\frac{9}{2^{k/2}},
\end{equation}
where we used the fact that $n+m-4>k/2$. Now we use for the third time Theorem \ref{Matveev} with $t:=3$ and the parameters
\[
(\eta_{1},b_{1}):=(\g,\ell), \quad (\eta_{2},b_{2}): =(2,-(n+m-3)) \quad \text{and} \quad (\eta_{3},b_{3}):=(\sqrt{2},-1).
\]
Note that $\eta_{1}$, $\eta_{2}$ and $\eta_{3}$ belong to the field $\mathbb{K}:=\mathbb{Q}(\sqrt{2})$ so that $D:=[\mathbb{K}:\mathbb{Q}]=2$. In addition, $h(\eta_{1})=(\log\g)/2$, $h(\eta_{2})=\log 2$ and $h(\eta_{3})=(\log 2)/2$, so we can choose $A_{1}:=\log\g$, $A_{2}:=2\log 2$ and $A_{3}:=\log 2$. Since $\ell<2n$ by \eqref{des-var}, we can take $B:=2n$ as we did before. It remains to prove that $\Lambda_{3}\neq 0$. If it were, then $\g^{2\ell}=2^{2n+2m-5}\in \Q$ which is not possible. Thus, $\Lambda_{3}\neq 0$.  It follows from Theorem \ref{Matveev} that
\begin{equation} \label{E4}
\log |\Lambda_{3}|>-1.1\times 10^{12}\log n,
\end{equation}
where we used the inequality $1+\log(2n)\leq 1.3 \log n$ for all $n\geq 420$. Comparing the above inequality \eqref{E4} with \eqref{forma3} leads to
\begin{equation}\label{aux1lamda3}
 k<3.2\times 10^{12}\log n.
\end{equation}
However, by Lemma \ref{n-in-terms-ofk} we have that
\[
\log n<\log(8.3\times 10^{30}k^{8}\log^{5}k)<71.2+8\log k+5\log\log k<22\log k
\]
for all $k>420$, which combined with \eqref{aux1lamda3} implies that
\[
k<7.1\times 10^{13}\log k.
\]
Combining the last inequality with Lemma \ref{n-in-terms-ofk} and \eqref{des-var}, we obtain the following absolute upper bounds:
\begin{equation}\label{cotasabsolutas}
k<2.6\times 10^{15}, \quad n<9.8\times 10^{161}  \quad \text{and} \quad \ell <2\times 10^{162}. 
\end{equation}
The next step involves reducing the above upper bounds by applying Lemma  \ref{reduce} once again. To this end, define 
\[
z_{3}:=\ell\log\g-(n+m-3)\log 2-\log\sqrt{2}
\]
and observe that \eqref{forma3} can be written as
\[
\left|e^{z_{3}}-1\right|<\frac{9}{2^{k/2}}.
\]
Note that $z_{3}\neq 0$. If $z_{3} > 0$, then by Lemma \ref{exponential}$(a)$ it follows that $|z_3| < | e^{z_3} -1 | < 9/2^{k/2}$.  If, on the contrary, $z_3 < 0$, one has $|e^{z_3} -1| <9/2^{k/2}<1/2$ since $k\ge 420$. Hence, from Lemma \ref{exponential}$(b)$, it holds that $|z_3| < 2| e^{z_3} -1 | <18/2^{k/2}$. In any case, the inequality $|z_3| <18/2^{k/2}$ is valid.  Replacing $z_3$ in the aforementioned inequality by its formula and dividing it across by $\log 2$, one obtains that 
\begin{equation}\label{red3}
0<\left|\ell\left(\frac{\log\g}{\log 2}\right)-(n+m-3)-\frac{1}{2}\right|<26\cdot (\sqrt{2})^{-k}.
\end{equation}
We now apply Lemma \ref{reduce} to inequality \eqref{red3} with the choices
\[
\tau:=\frac{\log\g}{\log 2}, \quad \mu:=-\frac{1}{2}, \quad A:=26 \quad \text{and} \quad B:=\sqrt{2}.
\]
It is clear that $\tau$ is an irrational number. Note that $\ell<M:=2\times 10^{162}$ by \eqref{cotasabsolutas}. By using \emph{Mathematica}, we found that $p_{302}/q_{302}$ is a convergent of $\tau$ with $q_{302}>6M$. A direct computation yields $\epsilon:=\left|\mu q_{302}\right\|-M\left\|\tau q_{302}\right\|>0.499$. It then follows by Lemma \ref{reduce} that
\[
k<\frac{\log \left(Aq_{303}/\epsilon\right)}{\log B}\leq 1110.
\]
Thus, $n<3.3\times 10^{59}$ and $\ell<2n<6.6\times 10^{59}$ by Lemma \ref{n-in-terms-ofk} and \eqref{des-var}.  The reduction process is then repeated using this refined upper bound for $\ell$. Indeed, applying Lemma \ref{reduce} to the inequality \eqref{red3} with $M:=6.6\times 10^{59}$ as the new upper bound for $\ell$ yields $k\leq 420$, which contradicts our assumption that $k > 420$. This completes the
analysis of the case when $k$ is large and therefore the proof of Theorem \ref{Main_Theorem}. 

\medskip

\textbf{Acknowledgements}.  J.~J.~B.~was supported in part by Project VRI ID 6115 (Universidad del Cauca). P.~D.~wants to thank the Universidad del Cauca in Popay\'an for hosting and his visit was funded by \emph{``Vicerrectoria de Investigaciones de la Universidad del Cauca"}. J.~L.~H.~was partially supported in part by Project VRI ID 6253 (Universidad del Cauca).



\begin{thebibliography}{100}

\bibitem{A1} M.~A.~Alekseyev, \emph{On the intersections of Fibonacci, Pell, and Lucas numbers}, Electronic Journal of Combinatorial Number Theory, Integers {\bf 11} (2011), 239--259.



\bibitem{BG2016} J.~J.~Bravo and C.~A. ~G\'omez, \emph{Mersenne $k$-Fibonacci numbers},  Glas.~Mat.~Ser.~III, {\bf 51} (2016), no. ~2,  307--319.

\bibitem{BGH-19} J. J. Bravo, C. A. G\'omez and J. L. Herrera, \emph{On the intersection of $k$-Fibonacci and Pell numbers}, Bull. Korean Math. Soc. {\bf  56} (2019), no. 2, 535--547.

\bibitem{BGL2016} J.~J.~Bravo, C.~A.~G\'omez  and F.~Luca, \emph{Powers of two as sums of two $k$-Fibonacci numbers}, Miskolc Math. Notes, {\bf 17} (2016), no.~1, 85--100.

\bibitem{BGL1} J. J. Bravo, C.~A. ~G\'omez  and F. Luca, \emph{A Diophantine equation in $k$-Fibonacci numbers and repdigits}, Colloq. Math. {\bf 152} (2018), no. 2, 299--315.

\bibitem{BHL_2021_JNT} J.~J.~Bravo, J.~L.~Herrera and F.~Luca, \emph{Common values of generalized Fibonacci and Pell sequences}, J. Number Theory, {\bf 226} (2021), 51--71.

\bibitem{BL1} J.~J.~ Bravo and F. ~Luca, \emph{On a conjecture about repdigits in $k$-generalized Fibonacci sequences}, Publ. Math. Debrecen, {\bf 82} (2013), no.~3--4, 623--639.

\bibitem{BravoLuca2015} J.~J.~Bravo and F.~Luca, \emph{Repdigits as sums of two  $k$-Fibonacci numbers}, Monatsh. Math., {\bf 176} (2015), no. 1, 31--51.















\bibitem{Cohn96} J. H. E. Cohn, \emph{Perfect Pell powers}, Glasgow Math.  J. {\bf 38} (1996), no. 1, 19--20.

\bibitem{Dasdemir} A. Da\c{s}demir and M. Varol, \emph{Lucas numbers which are products of two Pell numbers}, Fibonacci Quart., {\bf 63}(1) (2025), 78-83.

\bibitem{DdamuLucaRako2016} M. Ddamulira, F. Luca and M. Rakotomalala, \emph{Fibonacci numbers which are products of two Pell numbers}, Fibonacci Quart., {\bf 54}, no. 1, (2016), 11--18.

\bibitem{DD} G.~P.~Dresden and Z.~Du, \emph{A simplified Binet formula for $k$-generalized Fibonacci numbers}, J. Integer Seq., {\bf 17} (2014), no. 4, Article 14.4.7, 9 pp.

\bibitem{DP} A. Dujella and A. Peth\H o, \emph{A generalization of a theorem of Baker and Davenport}, Quart. J. Math. Oxford Ser., \textbf{49} (1998), no.~195, 291--306.




\bibitem{GL} S.~Guzm\'an and F.~Luca, \emph{Linear combinations of factorials and $S$-units in a binary recurrence sequence}, Ann. Math. Qu\'e., {\bf 38} (2014), no.~2, 169--188.

\bibitem{HY} L. K. Hua and Y. Wang, \emph{Applications of number theory to numerical analysis\/}, Translated from Chinese. Springer-Verlag, Berlin-New York; Kexue Chubanshe (Science Press), Beijing, 1981.

\bibitem{Kos01} T. Koshy, \emph{Fibonacci and Lucas Numbers with Applications, Volume 2}, Pure and Applied Mathematics, Wiley--Interscience Publications, New York, 2019.

\bibitem{Kos14} T. Koshy, \emph{Pell and Pell-Lucas Numbers with Applications}, Springer, New York, 2019.





\bibitem {Ljunggren42} W. Ljunggren, \emph{Zur Theorie der Gleichung $x^2+1=Dy^4$}, Avh.  Norske Vid Akad. Oslo {\bf 5} (1942).

\bibitem{Matveev} E. M. Matveev, \emph{An explicit lower bound for a homogeneous rational linear form in the logarithms of algebraic numbers}, II, Izv. Ross. Akad. Nauk Ser. Mat., \textbf{64} (2000), no.~6, 125--180; translation in Izv. Math., \textbf{64} (2000), no.~6, 1217--1269.

\bibitem{Mi} E. ~P. ~Miles, ~Jr., \emph{Generalized Fibonacci numbers and associated matrices}, Amer. Math. Monthly, {\bf 67} (1960), no. 8, 745--752.

\bibitem{Mil} M.~ D. ~Miller, \emph{Mathematical Notes: On Generalized Fibonacci Numbers}, Amer. Math. Monthly, {\bf 78} (1971), no. 10, 1108--1109.

\bibitem{Petho91} A. Peth\H{o}, \emph{The Pell sequence contains only trivial perfect powers}. Sets, graphs and numbers (Budapest, 1991), 561--568, Colloq. Math. Soc. J\'anos Bolyai, 60, North--Holland, Amsterdam, 1992. MR94e:11031.

\bibitem{Rihane1}S. Rihane, A. Youssouf, and A. E. Habibi. \emph{Fibonacci numbers which are products of three Pell numbers and Pell numbers which are products of three Fibonacci numbers}, Bol. Soc. Mat. Mex., {\bf 26} (2020), no. 3, 895-910.

\bibitem{Rihane2}S. Rihane, \emph{On k-Fibonacci and k-Lucas numbers written as a product of two Pell numbers}, Bol. Soc. Mat. Mex., {\bf 30} (2024), no. 1, 20.

\bibitem{Sloane} N.~Sloane, \emph{The On-line Encyclopedia of Integer Sequences} (OEIS), OEIS Foundation, Incorporated, 1996.

\bibitem{Waldshmidt} M.~Waldshmidt, \emph{Diophantine Approximation on Linear Algebraic Groups: Transcendence Properties of the Exponential Function in Several Variables}, Springer-Verlag Berlin Heidelberg, Berlin, 2000.

\bibitem{DAW} D. A. Wolfram, \emph{Solving generalized Fibonacci recurrences}, Fibonacci Quart., \textbf{36} (1998), no.~2, 129--145.

\end{thebibliography}
\end{document}